\documentclass{amsart}
\usepackage{latexsym,amssymb,amsfonts,amsmath, graphicx, color,enumerate}
\usepackage{dsfont}
\usepackage[latin1]{inputenc}
\usepackage{hyperref}
\numberwithin{equation}{section}
\newtheorem{theo}{Theorem}

\newtheorem{lemma}[theo]{Lemma}
\newtheorem{prop}[theo]{Proposition}
\newtheorem{cor}[theo]{Corollary}
\newtheorem{defi}[theo]{Definition}

\theoremstyle{definition}

\newtheorem{rem}[theo]{Remark}

\newtheorem{exa}[theo]{Example}

\numberwithin{theo}{section}
\DeclareMathOperator*{\esssup}{ess\,sup}

\title{Parabolic systems with coupled boundary conditions}
\author{Stefano Cardanobile}
\address{Bernstein Center for Computational Neuroscience,  Hansastra{\ss}e 9A, D-79104 Freiburg, Germany}
\email{stefano.cardanobile@bccn.uni-freiburg.de}
\author{Delio Mugnolo}
\address{Institut f\"ur Analysis, Universit\"at Ulm, Helmholtzstra{\ss}e 18, D-89081 Ulm, Germany}
\email{delio.mugnolo@uni-ulm.de}
\thanks{The authors would like to thank Robin Nittka and Wolfgang Arendt (Ulm) for interesting discussions.}

\begin{document}

\begin{abstract}
We consider elliptic operators with operator-valued coefficients and discuss the associated parabolic problems. The unknowns are functions with values in a Hilbert space $W$. The system is equipped with a general class of coupled boundary conditions of the form $f_{|\partial\Omega}\in \mathcal Y$ and $\frac{\partial f}{\partial \nu}\in {\mathcal Y}^\perp$, where $\mathcal Y$ is a closed subspace of $L^2(\partial\Omega;W)$. We discuss well-posedness and further qualitative properties, systematically reducing features of the parabolic system to operator-theoretical properties of the orthogonal projection onto $\mathcal Y$.
\end{abstract}

\maketitle

\section{Introduction}

In this paper we investigate diffusive equations in a vector valued setting. More precisely, we discuss a problem of the form
$$\frac{\partial u}{\partial t}(t,x)={\mathcal L} u(t,x),\qquad t\geq 0,\;x\in\Omega\subset{\mathbb R}^n, $$
where the unknown $u$ takes values in a separable Hilbert space $W$ and $\mathcal L$ is an elliptic operator with operator-valued symbol.

Operators on vector-valued spaces have been often considered in recent years. We mention the important investigations performed in this field by Amann in~\cite{Ama84,Ama01} and by Denk--Hieber--Pr\"uss in~\cite{DenHiePru03,DenHiePru07}, cf.\ also references therein. They consider general differential (or pseudo-differential) operators of order $2m$ on vector-valued spaces, whereas for the sake of simplicity we restrict ourselves to the elliptic case with $m=1$ only in the present article.

Operators with operator-valued symbol also appear in theoretical physics, in particular in quantum mechanics and molecular dynamics. 
There, the relevant state space is often of the form $$L^2({\mathbb R}^{3l})\otimes L^2({\mathbb R}^{3k})\equiv L^2({\mathbb R}^{3l},L^2({\mathbb R}^{3k})),$$
where $L^2({\mathbb R}^{3l})$ and $L^2({\mathbb R}^{3k})$ are the configurations spaces of different parts of the system, e.g.\ in the time-dependent Born--Oppenheimer theory.

We also add that a variational approach has been proposed in~\cite{CarMug08}, where differential operators (with values in finite-dimensional spaces) are rather described by means of a matrix formalism.

The first systematic investigation of elliptic systems with possibly coupled boundary conditions probably goes back to the classic work by Agmon--Douglis--Nirenberg. In~\cite{AgmDouNir64} they have considered quite general boundary conditions that, in our context of second-order operators, essentially read
\begin{equation}\label{adn}
Af_{|\partial\Omega}+B\frac{\partial f}{\partial \nu}=0 
\end{equation}
where (up to slight generalizations) $A,B$ are differential operators from $L^2(\Omega;W)$ to $L^2(\partial\Omega;W)$ that satisfy an involved ``complementing'' condition, cf.~\cite[\S~I.2]{AgmDouNir64}. Successively, systems with coupled boundary conditions have been only occasionally investigated, mostly in the case of reaction-diffusion equations or of Sturm--Liouville problems (so that $W={\mathbb C}^n$). They usually feature some kind of generalized, coupled Robin conditions, see e.g.~\cite{Leu82,Ama86,Ama88}, or~\cite[Chapt.~13]{Zet05}, respectively. Another concrete class of systems that fit the setting of Agmon--Douglis--Nirenberg has been thoroughly investigated by Ali Mehmeti and Nicaise in the context of the theories of \emph{interaction problems} and \emph{elliptic operators on polygonal domains}, cf.~\cite{Nic93,AliNic93,Ali94} and references therein -- see also~\cite[\S~III.4.4]{Sho94} for a discussion of the associated elliptic problem. In the 1-dimensional case (i.e., $\Omega\subset \mathbb R$), mathematical investigations of systems with coupled boundary conditions (i.e., \emph{networks}) go back to Lumer~\cite{Lum80}. The field of elliptic operators on 1-dimensional networks has become so popular over the last 10 years, both in applied analyis and mathematical physics, that an overview of the existing literature is impossible in this limited introduction. Such an overwhelming work has been successfully performed in~\cite{Kuc08}, though.

In this paper we discuss a setting for mixing boundary conditions that allows for truly coupled dynamics of the system, even when the coupling is only defined on the boundary. In the prototypical case, we discuss boundary conditions that can be formulated as
\begin{equation}\label{kuch}
f(z)\in Y\qquad\hbox{and}\qquad \frac{\partial f}{\partial \nu}(z)+Sf(z)\in Y^\perp\qquad\hbox{for a.e.\ }z\in\partial \Omega,
\end{equation}
where $Y$ is any subspace of $W$, or more generally as
\begin{equation*}
f_{|\partial\Omega}\in {\mathcal Y}\qquad\hbox{and}\qquad \frac{\partial f}{\partial \nu}+Rf_{|\partial\Omega}\in {\mathcal Y}^\perp,
\end{equation*}
where $\mathcal Y$ is any subspace of $L^2(\partial\Omega;W)$. Here $R$ and $S$ are bounded linear operators on $L^2(\partial\Omega;W)$ and on $W$, respectively. In the case of $1$-dimensional domains $\Omega$ and finite dimensional $W$, essentially the same boundary conditions have been thoroughly discussed in~\cite{Kuc04,FulKucWil07} in the context of quantum graphs. Kuchment has also proved that if $A,B$ define via~\eqref{adn} a self-adjoint elliptic operator, then the same boundary conditions can also be formulated as in~\eqref{kuch}, for suitable $Y$ and $R$, cf.~\cite[Thm.~6]{Kuc04}. His proof carries over to the case of higher dimensional domains $\Omega$ if still ${\rm dim}\, W<\infty$.

As a matter of fact, our  boundary conditions form a proper subset of the class considered in~\cite{AgmDouNir64} for the case ${\rm dim}\, W<\infty$ and extended to the infinite dimensional vector-valued case in ~\cite{DenHiePru03}. Still, our formulation seems to have some advantages over the general Agmon--Douglis--Nirenberg's one. In particular, qualitative properties can be simply characterized by linear algebraic (resp., operator theoretical) properties (if $W$ is finite or resp. infinite dimensional) of the \emph{lower dimensional} orthogonal projections of $W$ onto $Y$. 

To do so, we pursue in this paper a variational approach that is similar to that applied in~\cite{MerNic98,CarMugNit08} to the strongly coupled, 1-dimensional case of $\Omega=(0,1)$. Parabolic systems in variational form have been discussed at an abstract level also in~\cite{Ouh99}.
In Section~\ref{general} we prove a basic well-posedness result by discussing properties of a suitable sesquilinear form and identifiy the associated operator $A$. Describing the behavior of the semigroup generated by $A$ by means of relevant invariance properties is the main goal of the subsequent sections. Smoothing properties of this semigroups are discussed in Section~\ref{regularity}.

In Section~\ref{symmetries} we introduce a class of subspaces of the state space $L^2(\Omega;W)$ defined pointwise by a linear algebraic relation. Following~\cite{CarMugNit08}, their invariance under the diffusion semigroup can be interpreted as a symmetry property of the associated system of Schr\"odinger equations.
Section~\ref{posdom} is mostly devoted to positivity and $L^\infty$-contractivity, and more generally to the invariance of sets of functions pointwise belonging to an order interval of $W$.
In Section~\ref{dominff} we discuss domination issues, showing that the evolution of the diffusion problem with Dirichlet boundary conditions is minimal in the class considered throughout the paper; and that the evolution of the diffusion problem with Neumann boundary conditions is maximal among those problems whose boundary space $Y$ is a closed ideal of $W$. {In~\cite{AreWar03} Arendt and Warma have showed that, in the scalar valued setting, the only boundary conditions that a) are local and b) give rise to a semigroup dominated by the Neumann semigroup \emph{and} dominating the Dirichlet semigroup are of Robin type. Similarly, we can prove that the only boundary conditions that give rise to a semigroup dominated by the Neumann semigroup \emph{and} dominating the Dirichlet semigroup in the vector-valued case are necessarily of \emph{decoupled} Robin type.}

%In Section~6 we discuss the issue of irreducibility, i.e., invariance of closed subspaces, and show that this differs from special cases met in applications. In the Appendix, we formulate the general principle that allows us to extend known results to the vector-valued case and prove a technical result.

We will use the notion of Hilbert lattice over and over and refer to~\cite[Chapter~C]{Nag86} and~\cite[Chapter~2]{Mey91} for a general treatment of this theory.

\begin{rem}
After the first draft of this article was completed, we discovered at the Conference on ``Semigroups and Evolution Equations'' held in T\"ubingen in November 2008 that comparable results have been independently obtained in a manuscript by Ulrike Kant, Tobias Klau{\ss}, J\"urgen Voigt, and Matthias Weber~\cite{KanKlaVoi09}. However, their focus is on compact quantum graphs, i.e., on finite $1$-dimensional networks, whereas we rather consider a setting based on general Hilbert-space-valued function spaces.
\end{rem}

\section{General setting}\label{general}

Throughout this paper, $\Omega\subset\mathbb R^N$ is a domain with Lipschitz boundary and $W$ a separable Hilbert space. Define $H := L^2(\Omega;W)$ and $V := H^1(\Omega;W)$% as well as $\partial V := H^\frac{1}{2}(\partial\Omega;W)$ and $\partial H := L^{2}(\partial\Omega;W)$
: these are apparently Hilbert spaces. Observe that, due to separability of $W$, these spaces can be defined in the usual way by means of Bochner integrals, see e.g.~\cite{DieUhl77}. More precisely, a function is in $H$ if and only if it is weakly measurable and the scalar-valued function $\Omega\ni x\mapsto \|f(x)\|_W\in\mathbb R$ is of class $L^2(\Omega;{\mathbb R})$.

In fact, most common results concerning boundary regularity and Gau{\ss}-Green formulae for scalar functions remain valid for vector-valued functions, as the usual proofs can be extended to the vector-valued case by the Pettis' measurability theorem.
Although such results are known, we discuss in the Appendix (Section~\ref{extensapp}) an abstract device that permits to extend them to the vector-valued case in a systematic way.

If we consider a closed linear subspace $\mathcal Y \subset {L^2(\partial\Omega;W)}$, by Example~\ref{tracexa} it is possible to define
\begin{equation}\label{VE}
V_{\mathcal Y}:=\{f\in V:f_{| \partial\Omega}\in \mathcal Y\},
\end{equation}
which is a closed subspace of $V$, due to the continuity of the trace operator. By $f_{|\partial\Omega}$ we denote the trace of $f$ defined in accordance with Example~\ref{tracexa}.

Several common examples of closed subspaces of $V$
appearing in applications are in fact of
the above form, most notably $V_{{L^2(\partial\Omega;W)}} = V$ and
$V_{\{0\}} = H^{1}_0(\Omega;W)$. If $\Omega\subset\mathbb R$ is a bounded interval, then also periodic and anti-periodic boundary conditions fit into our framework.

\begin{rem}
For a separable Banach space $E$ we denote by
$L^\infty(\Omega; {\mathcal L}_s E)$ the space of all functions $f: \Omega \to {\mathcal L}(E)$ such that $\omega \mapsto f(\omega)x$ is of class $L^\infty(\Omega,E)$ for every $x \in E$. This notation has been introduced in~\cite{AreTho05}, where a thorough investigation of vector-valued $L^p$-spaces has been performed.
\end{rem}

Fix $D\in L^\infty(\Omega;{\mathcal L}_s(W^n))$ and $R \in{\mathcal L}(V_{\mathcal Y};{L^2(\partial\Omega;W)})$ and consider the sesquilinear form
$a:=a_{\mathcal Y}:V_{\mathcal Y}\times V_{\mathcal Y}\to\mathbb C$ defined by
\begin{eqnarray}\label{aform}
a(f,g)&:=&\int_\Omega \left(D(x)\nabla f(x) | \nabla g(x)\right)_{W^n} dx+\int_{\partial\Omega} \left((R f)(z) | g(z)\right)_W d\sigma(z) \nonumber\\
&\equiv &\sum_{h,k=1}^n \int_\Omega \left(D_{hk}(x) \partial_k f(x)|{\partial_h g(x)}\right)_W dx+\int_{\partial\Omega} \left((R f)(z) | g(z)\right)_W d\sigma(z). \nonumber
\end{eqnarray}

\begin{rem}
While networks seem to be a simple application of our theory, properties of vector-valued diffusion problems do sometimes differ from those of network parabolic problems. If e.g.\ $W=\ell^2(E)$, where $E$ is an infinite countable set, then it is clear that $L^2(0,1;W)$ is isometrically isomorphic to $L^2(\mathbb R;{\mathbb C})$. 
Still, two given functions need not have disjoint support in $L^2(0,1;W)$ when their versions in $L^2(\mathbb R;{\mathbb C})$ do. 
In particular, the form $a_{\mathcal Y}$ is local for all $\mathcal Y \subset {L^2(\partial\Omega;W)}$, regardless of the locality of the operators $D_{hk}(x)$, $x\in\Omega$. This should be compared with the setting discussed in~\cite{CarMugNit08} for network equations.

% 
% The form $a$ is in general not local, i.e., it may happen that $a(f,g)\not=0$ for some $f,g\in V_{\mathcal Y}$ having disjoint support, even in the simpler case of $R=0$. In this case, $a$ is given by
% $$a(f,g)=\sum_{h,k=1}^n \int_\Omega \left(D_{hk}(x) \partial_k f(x)| {\partial_h f(x)}\right)_{W^n} dx.$$
% \textcolor{red}{If $W=L^2(X)$, then a sufficient and necessary condition for locality of $a$ is that $D_{hk}(x)\in {\mathcal L}(W)$ is a multiplication operator for all $h,k=1,\ldots,n$ and all $x \in \Omega$, in analogy with the case of strongly coupled diffusion system, see~\cite{CarMugNit08}.}
\end{rem}

Throughout the paper we impose the assumption that the ellipticity condition
$${\rm Re}(D(x)\xi| \xi)_{W^n}\geq \gamma \| \xi\|_{W^n}^2\qquad\hbox{ for all } \xi\in W^n\hbox{ and a.e.\ }x\in\Omega$$	
is satisfied for some $\gamma>0$.

\begin{prop}\label{form}
Let $\mathcal Y$ be a closed subspace of $L^2(\Omega;W)$. Let $R\in{\mathcal L}(V_{\mathcal Y},{L^2(\partial\Omega;W)})$. Then $a$ is $H$-elliptic, continuous, and densely defined.

Furthermore, the form is accretive (resp., coercive) if 
\begin{equation}\label{characcr}
{\rm Re}(Rf | \hat{T}f)_{{L^2(\partial\Omega;W)}}\geq 0 \hbox{ (resp., }>q\|f\|_{{L^2(\partial\Omega;W)}}\hbox{ for some }q>0)\qquad \hbox{ for all }f\in V_{\mathcal Y}.
\end{equation}
It is symmetric if and only if $D(x)$ is self-adjoint for a.e.\ $x\in\Omega$ and moreover
\begin{equation}\label{charsymm}
(Rf |  \hat{T} g)_{{L^2(\partial\Omega;W)}}=\overline{( R g | \hat{T} f)}_{{L^2(\partial\Omega;W)}}\qquad\hbox{ for all }f,g\in V_{\mathcal Y},
\end{equation}
where $\hat{T}$ denotes the vector-valued trace operator defined in Example~\ref{tracexa}.
\end{prop}

\begin{proof}
Boundedness of the trace operator $\hat{T}:V\to{L^2(\partial\Omega;W)}$ and the Cauchy--Schwarz inequality yield
$$|a(f,g)| \leq \esssup_{x\in\Omega}\| D(x)\|_{{\mathcal L}(W^n)}\| \nabla f\| _H \| \nabla g\| _H
+M\| R \|_{\mathcal L(V_{\mathcal Y},{L^2(\partial\Omega;W)})}  \| f\| _{V}\| g\| _{V},$$
hence $a$ is continuous. It is $H$-elliptic, since $a$ can be considered as the sum
$$a(f,g)=a_1(f,g)+a_2(f,g):=\int_\Omega \left(D(x)\nabla f(x) | \nabla g(x)\right)_{W^n} dx
+\int_{\partial\Omega} \left((R f)(z) | g(z)\right)_W d\sigma(z),$$
where by assumption $a_1$ is $H$-elliptic and $a_2$ is bounded on, say, $V\times H^{\frac{3}{4}}(\Omega;W)$. Thus, by~\cite[Lemma~2.1]{Mug08} also the perturbed form $a$ is $H$-elliptic. 
The remaining assertions are clear due to accretivity and symmetry of the leading term $a_1$.
\end{proof}

Thus, the operator $A_{\mathcal Y}$ associated with $(a_{\mathcal Y},V_{\mathcal Y})$ generates an analytic semigroup $(e^{tA_{\mathcal Y}})_{t\geq 0}$ on $H$. 

This semigroup is self-adjoint if and only if $D(x)$ is self-adjoint for a.e.\ $x\in\Omega$ and~\eqref{charsymm} holds. It is contractive (resp., uniformly exponentially stable) if~\eqref{characcr} holds.

Moreover, provided that $\Omega$ is bounded, by Aubin's Lemma (see~\cite{Aub63}) the semigroup is compact if and only if $W$ is finite dimensional. If this is the case,  and if $R=0$, then $(e^{tA_{\mathcal Y}})_{t\geq 0}$ is uniformly exponentially stable if and only if $A_{\mathcal Y}$ is invertible, i.e., if and only if $a(f,f)=0$ implies $f=0$. Under our ellipticity assumption, this is possible if and only if $\nabla f=0$ implies $f=0$. Clearly, this is only possible if the boundary conditions force each $f\in V_{\mathcal Y}$ to vanish on some non-null subset of the boundary, i.e., if and only if 
\begin{equation}\label{conditstab}
{\mathcal Y}=\{0\}\oplus L^2(\partial\omega;W)\qquad \hbox{for some }\partial\omega\subset\partial\Omega, 
\end{equation}
with $\partial\Omega\setminus\partial\omega$ of non-zero measure. Observe that if we drop the compactness assumption
but~\eqref{conditstab} holds, then by~\cite[Exa.~V.2.23]{EngNag00} the semigroup is still \emph{strongly} stable.

%\footnote{\textcolor{red}{in realt\`a mi sono reso conto che aubin prova solo un'implicazione (se $W\hookrightarrow Z$ \`e compatta, allora $H^1(\Omega,W)\hookrightarrow L^2(\Omega,Z)$ \`e compatta). vale anche il viceversa, nei casi non banali (tipo $\Omega=\emptyset$)? a me sembra ragionevole, ma non riesco a provarlo}
%\textcolor{blue}{Mi pare il viceversa vale per qualsiasi $\Omega$ non banale. Tentativo di dimostrazione: dobbiamo provare che $H^1 \hookrightarrow L^2$ compatta implica $W$ finito-dimensionale.
%Si fissi una successione di norma unitaria $w_n \subset W$ e una funzione $\phi \in H^1(\Omega, \mathbb C)$, anche essa di norma unitaria. 
%Evidentemente la successione $w_n \otimes \phi \subset H^1(\Omega, W)$ ha norma unitaria e converge in $H^1$ o in $L^2$ se e solo se $w_n$ converge in $W$.
%Dato che l'immersione \'e compatta, $w_n \otimes \phi $ ha una sottosuccessione convergente in $L^2$, quindi $w_n$ ha una sottosuccessione convergente. Per l'arbitrariet\'a di $w_n$ concludiamo che $W$ ha dimensione finita}}. 

\begin{rem}
It follows by standard perturbation results, see e.g.~\cite[Lemma~2.1]{Mug08}, that similar results also hold by considering lower order terms. Let e.g.\ $n\ge 2$ and $p\ge\frac{2n}{n-2}$, so that $H^1(\Omega)\hookrightarrow L^p(\Omega)$. It is known that the same embedding is satisfied by vector-valued Sobolev spaces. Let additionally $C \in L^\frac{p}{p-2}(\Omega;{\mathcal L}(W))$. Then by H\"older's inequality also the form defined by
$$\int_\Omega \left(D(x)\nabla f(x) | \nabla g(x)\right)_{W^n} dx+\int_\Omega (C(x)f(x)|g(x))_W dx+\int_{\partial\Omega} \left((R f)(z) | g(z)\right)_W d\sigma(z)$$
for $f,g\in V_{\mathcal Y}$ is elliptic, continuous and densely defined, and we conclude that the associated diffusion problem is well-posed. Weakly coupled systems modelled by diffusion on vector-valued spaces have been considered by several authors, for instance in order to discuss molecular motors -- see e.g.\ \cite{ChiHilKin08} and references therein. 
Relevant properties of molecular models of this kind, like $L^1$-contractivity or positivity, will be discussed in Section~\ref{posdom}.
\end{rem}

Since $D(A_{\mathcal Y})\subset V_{\mathcal Y}$, all functions $f \in D(A_{\mathcal Y})$
satisfy the boundary condition $f_{|\partial \Omega} \in {\mathcal Y}$. With only this boundary condition, though, the system is in general underdetermined. In fact, the following holds.

\begin{prop}
Assume $D\in C^1(\overline{\Omega};{\mathcal L}(W^n))$ and the domain $\Omega$ to have $C^1$-boundary. Under the assumptions of Proposition~\ref{form} let
\begin{eqnarray*}
D(B_{\mathcal Y})&:=&\left\{f\in C^2(\overline{\Omega};W): f_{| \partial\Omega}\in {\mathcal Y}\mbox{ and }\frac{\partial_D f}{\partial\nu}+Rf\in {\mathcal Y}^\perp\right\},\\
B_{\mathcal Y}f&:=&\nabla \cdot(D\nabla f),
\end{eqnarray*}
where $\nu(z)$ denotes the outer normal at $z\in\partial\Omega$ and 
$$\frac{\partial_D f}{\partial\nu}(z):=\left(D(z)\nabla f(z)|\nu(z)\right)_{{\mathbb C}^n}:= 
\sum_{i,j=1}^n \nu_i(z) D_{ij}(z){\partial_j } f(z)$$
is the (vector-valued) conormal derivative w.r.t.\ $D$ at $z\in\partial\Omega$. Then $B_{\mathcal Y}\subset A_{\mathcal Y}$. On the other hand, if $f\in D(A_{\mathcal Y})\cap C^2(\overline{\Omega};W)$ then
$$\frac{\partial_D f}{\partial\nu}+R f\in\mathcal Y^\perp.$$
\end{prop}

\begin{proof}
To begin with, observe that if $f\in C^2(\overline{\Omega};W)$, then it belongs to $H^2(\Omega;W)$ and hence by the general theory of Sobolev spaces and by Theorem~\ref{extension} it admits well-defined (vector-valued) trace and normal derivative of class $L^2(\partial\Omega;W)$. In particular, Gau{\ss}--Green's formulae hold.

The proof mimicks~\cite[Rem.~3.1.6]{Are06}. Observe that by Theorem~\ref{extension} we can replicate the proof of~\cite[Thm.~IX.8]{Bre83} and conclude that vector-valued test functions are dense in $H^1(\Omega;W)$.

If $f \in D(B_{\mathcal Y})$ and $h\in V_{\mathcal Y}$, then by the Gau{\ss}-Green formulae
\begin{eqnarray*} -(B_{\mathcal Y} f | h)_H &=&\int_\Omega \left(D(x)\nabla f(x) | \nabla h(x)\right)_{W^n} dx-\int_{\partial\Omega}\left(\frac{\partial_D f}{\partial\nu}| h(z) \right)_W d\sigma(z) \\
&=&\int_\Omega \left(D(x)\nabla f(x) | \nabla h(x)\right)_{W^n} dx+\int_{\partial\Omega}\left((R f)(z) | h(z)\right)_W d\sigma(z) =a(f,h),
\end{eqnarray*} 
hence $B_{\mathcal Y}\subset A_{\mathcal Y}$.

Let conversely $f\in D(A_{\mathcal Y})\cap C^2(\overline{\Omega};W)$ and $h \in V_{\mathcal Y}$. Then by definition of $A_{\mathcal Y}$
\begin{eqnarray*} -(A_{\mathcal Y} f | h)_H & =& a(f,h) \\
&=&\int_\Omega \left(D(x)\nabla f(x) | \nabla h(x)\right)_{W^n} dx+\int_{\partial\Omega}\left((R f)(z) | h(z)\right)_W d\sigma(z) \\
&=&\int_{\partial\Omega} \left(\frac{\partial_D f}{\partial\nu} | h(z)\right)_{W}dx - \int_\Omega \left(\nabla\cdot (D\nabla f) (x) | h(x)\right)_W dx\\
&&\qquad\qquad + \int_{\partial\Omega}\left((Rf)(z) | h(z)\right)_W d\sigma(z).
\end{eqnarray*} 
In particular, this computation holds for $h\in H^{1}_0 (\Omega; W)$. It follows that
$A_{\mathcal Y}f=\nabla\cdot(D\nabla f)$ for all $f \in D(A_{\mathcal Y})$. As a consequence, for all $h \in C^1(\overline{\Omega};{W})$ such that $h_{|\partial \Omega} \in {\mathcal Y}$ we obtain
$$\int_{\partial\Omega} \left(\left(\frac{\partial_D f}{\partial\nu}+R f\right)(z) | h(z)\right)_W d\sigma(z)=0.$$
This implies 
$$\frac{\partial_D f}{\partial \nu}+R f\in {\mathcal Y}^\perp$$
because the trace operator has dense range in ${L^2(\partial\Omega;W)}$.
\end{proof}

\section{Regularity properties}\label{regularity}

In the previous section we have showed generation properties for the operator $A_{\mathcal Y}$ associated with the form $a_{\mathcal Y}$, i.e., that $A_{\mathcal Y}$ generates an analytic semigroup $(e^{tA_{\mathcal Y}})_{t \geq 0}$. In this and the next sections we will discuss the issue of invariance of relevant subspaces or, more generally, subsets of $H$ under the action of $(e^{tA_{\mathcal Y}})_{t \geq 0}$. To do so, throughout the remainder of this article we specialize by considering a particular class of boundary conditions.

Fix a closed linear subspace $Y \subset W$ and consider the Hilbert subspace
\begin{equation}\label{formdomain}
V_Y:=\{f \in H^1(\Omega; W): f(z) \in Y\hbox{ for a.e.\ }z\in\partial\Omega\}\subset H^1(\Omega;W).
\end{equation}
(We emphasize the difference between this space and that introduced in~\eqref{VE}).

Diffusion equations on a domain with Dirichlet or Neumann boundary conditions, as well as on an ``open book'' with finitely many pages fit this setting.

\begin{exa}\label{orthogonaldirichlet}
Let $W = \mathbb C$, $D(x)=Id$ for a.e.\ $x\in\Omega$, and $R = 0$.
Then for $Y= \mathbb C$ one has $V_Y=H^1(\Omega;{\mathbb C})$ and
$$a(f,g)=\int_\Omega \nabla f(x)\overline{\nabla g(x)}dx,\qquad f,g\in H^1(\Omega;{\mathbb C}),$$ 
is the form associated with the Neumann Laplacian. If more generally $R f=\beta
f_{| \partial\Omega}$ for some $\beta:\partial\Omega\to\mathbb C$, then $a$ is
the form associated with the Robin Laplacian. Choosing $Y=\{0\}$ we obtain the form associated with the Dirichlet Laplacian.
\end{exa}

\begin{exa}\label{dihedron}
Let $\Omega$ be the open 2-dimensional halfplain $\{(x_1,x_2)\in {\mathbb R}^2:x_1\geq 0\}$, so that $\partial \Omega=\{0\}\times {\mathbb R}$. 
Let moreover $W={\mathbb C}^m$, $m\in\mathbb N$, and $Y$ be the subspace of $W$ spanned by $\mathbf 1=(1,\ldots,1)$. 
Then a function $f\in H^1(\Omega;{\mathbb C}^m)$ lies in $V_Y$ if and only if for all $z\in\partial \Omega$ there holds $f_1(z)=\ldots=f_m(z)$, i.e., if and only if $f$ is continuous along the ``binding'', giving rise to a Kirchhoff-type boundary condition. Such diffusion problems on open book structures have been discussed by several authors, see e.g.~\cite{Nic93,Ali94,FreWen04}. Although less physically motivated, one can also naturally consider the ``dual'' problem obtained replacing $Y$ by $Y^\perp$. Then the boundary condition along the binding is that $f_1(z)+\ldots+f_m(z)=0$. This condition is sometimes called \emph{anti-Kirchhoff} in the context of quantum graphs, cf.~\cite{Kuc04}.
\end{exa}

We first prove a regularity result for the semigroup generated by $A$ and that
the latter is given by a vector-valued integral kernel.

\begin{prop}\label{boundreg}
For $k\in\mathbb N$ assume $D\in C^k(\overline\Omega;{\mathcal L}(W^n))$ and
the domain $\Omega$ to have $C^{2k}$-boundary. Then the semigroup associated
with $a_Y$ maps $H=L^2(\Omega;W)$ into $H^{2k}(\overline{\Omega};W)$ for all $Y
\subset W$.
% \footnote{\textcolor{red}{probabilmente le ipotesi possono essere indebolite. per~\cite[cor. IX~13]{Bre83}, l'immersione $H^k(\Omega)\hookrightarrow L^\infty(\Omega)$ vale non appena $k>\frac{n}{2}$. occorre quindi risalire con teoremi di regolarit\`a al bordo semplicemente a $H^{k}(\Omega)$, con $\Omega\in C^k$ e $k>\frac{n}{2}$.}
% \textcolor{blue}{per me si pu\'o annerire, per\'o non so che dire sulla tua nota.}
% }
\end{prop}

\begin{proof}
It is known that an analytic semigroup on $H$ with generator $A$ maps $H$ into $\bigcap_{k=0}^\infty D(A^k)$. Thus, it suffices to show that $D(A^k)\subset H^{2k}(\Omega;W)$ for all $k\in\mathbb N$.

First, we prove that the assertion is true for $k=1$. In fact, if $u\in D(A)$, then $Au\in H$ and $u$ satisfies boundary conditions. 
Observe that if $\Omega$ is a half-plane, then the boundary conditions
\begin{equation*}
f(z)\in Y \quad\hbox{and}\quad \frac{\partial f}{\partial \nu}(z)+(Rf)(z)\in Y
\hbox{ for a.e.\ }z\in\partial\Omega
\end{equation*}
are still satisfied upon translating $u$ along the boundary. As a consequence, it is possible to mimic the proof of~\cite[Th\'eo.~IX.25]{Bre83} based on Nirenberg's technique of incremental quotients. 
(This also shows why this proof cannot be performed in the more general context of the boundary conditions considered in the previous section, which involved more general spaces $V_{\mathcal Y}$). Summing up, we obtain that $u\in H^2(\Omega;W)\cap V_Y$. 

Likewise, if $u\in D(A^k)$, then $Au\in D(A^{k-1})\subset H^{2(k-1)}(\Omega;W)$ by induction hypothesis.
By the second part of the assertion in~\cite[Th\'eo.~IX.25]{Bre83}, $u\in H^{2k}(\Omega;W)$. 
This concludes the proof in the case that $D_{hk}(x)=Id$ for a.e.\ $x\in \Omega$. The general case follows by the techniques presented e.g.\ in~\cite[\S~2.2.2]{Gri85} or~\cite[\S~5.6]{Eva98}.
\end{proof}

% \begin{rem}\label{domain}
% In fact, by the above mentioned boundary regularity results it is possible to show that 
% $$D(A)=\{f\in H^2(\Omega;W): f_{|\partial\Omega}\in Y \hbox{ and }\frac{\partial f}{\partial \nu}+Rf_{|\partial\Omega}\in Y^\perp\}$$ 
% and that $Af=\nabla \cdot (D\nabla f)$ for all $f\in D(A)$. This is already the case if we only assume $D\in C^1(\Omega;W^n)$ and $\Omega$ to have $C^1$-boundary. The proof is standard and we omit the details.
% \footnote{\textcolor{blue}{A me questo remark confonde solo le idee, lo togliamo?}}
% \end{rem}

% The following is a special case of a result that has been obtained in~\cite{MugNit08}, based on the classic Dunford--Pettis-criterion for representability of operators.
% 
% \begin{theo}\label{dplike}
% The following assertions are equivalent.
% \begin{enumerate}[(a)]
% \item The operator $T\colon L^1(\Omega;W) \to L^\infty(\Omega;W)$ is linear and bounded.
% \item There exists $k \in L^\infty(\Omega \times \Omega; \mathcal{L}_s(W))$ such that
% $$(Tf)(y) = \sigma\!\mbox{-}\!\int_{\Omega} k(x,y)f(x) d x$$
% for a.e. $y \in \Omega$, i.e., for every $v \in W$ and a.e. $y \in \Omega$ we have
% $$({(Tf)(y)}|v) = \int_{\Omega} (k(x,y)f(x)|v)_W d x.$$
% \end{enumerate}
% In this case, 
% $$\esssup_{(x,y)\in\Omega\times\Omega}\|k(x,y)\| = \esssup_{\genfrac{}{}{0pt}{}{y\in\Omega}{f\in L^1(\Omega; W)}} \frac{\|(Tf)(y)\|_{W}}{\|f\|}.$$
% \end{theo}
% 
% We remark that also the Hilbert--Schmidt property and, if $W$ is a Hilbert lattice, positivity and regularity can be characterized in a similar way, see~\cite[\S\S~5--6]{MugNit09}.

\begin{cor}
Assume $D\in C^j(\overline\Omega;{\mathcal L}(W^n))$ and the domain $\Omega$ to have $C^k$-boundary. If $2k>n$ and $2j\ge k$, then the semigroup associated with $a$ is given by an integral kernel.
\end{cor}

\begin{proof}
By Proposition~\ref{boundreg} and standard Sobolev embeddings (see e.g.~\cite[Cor.~IX.13]{Bre83}), the semigroup $(e^{tA_Y})_{t\ge 0}$ associated with $a$ maps $H$ into $H^k(\Omega;W)$ and hence into $L^\infty(\Omega;W)$. The same holds for the adjoint of $(e^{tA_Y})_{t\ge 0}$ and by duality $(e^{tA_Y})_{t\ge 0}$ maps $L^1(\Omega;W)$ into $L^\infty(\Omega;W)$. We conclude that the semigroup operators $e^{tA_Y}$ are given by an integral kernel of class $L^\infty(\Omega\times\Omega;{\mathcal L}_s(W))$ for all $t>0$ by~\cite[Thm.~3.3 and Prop.~4.3]{MugNit09}.
\end{proof}

\section{Symmetries}\label{symmetries}

After discussing smoothing properties of $(e^{tA_Y})_{t\ge 0}$ by means of invariance of dense subsets of $H$, we turn our attention to a different issue.

For an arbitrary closed convex subset $C_W\subset W$ we introduce closed convex subsets of $H=L^2(\Omega;W)$ and $\partial H=L^2(\partial\Omega;W)$ by
 $$ C_H:=\{f\in H: f(x)\in C_W \hbox{ for a.e.\ }x\in\Omega\}$$
and
 $$ C_{\partial H}:=\{f\in \partial H: f(z)\in C_W \hbox{ for a.e.\ }z\in\partial\Omega\}.$$
If $C_W$ is a subspace of $W$, then $C_H,C_{\partial H}$ are clearly subspaces, too. The connection between invariance of pointwise defined closed subspaces and the physical notion of symmetry has been discussed in~\cite[\S~5]{CarMugNit08} and in~\cite[\S~2.8]{Car08}.
. 
% \begin{rem}\label{projpa}

One sees that the orthogonal projection of $H=L^2(\Omega;W)$ onto ${C_H}$ is defined by
\begin{equation*}
(P_{{C_H}} f)(x) = P_{C_W}(f(x)), \qquad x \in \Omega,
\end{equation*}
where $P_{C_W}$ denotes the orthogonal projection of $W$ onto $C_W$.
% \end{rem}

For the sake of simplicity, in the following we assume the operator $R\in{\mathcal L}(V_Y,\partial H)$ to admit the factorization 
$$R= S\hat{T},$$
where $S\in{\mathcal L}(\partial H)$.

For the characterization of the invariance of a closed convex subset $C_H$ of $H$ it is convenient to use the criterion of E.M. Ouhabaz, cf.~\cite[Chapt.~2]{Ouh05} and its generalization proved in~\cite{ManVogVoi05}. It implies that $(e^{tA_Y})_{t\geq0}$ leaves $C_H$ invariant if and only if 
$$P_{{C_H}}V_Y \subset V_Y\qquad\hbox{and}\qquad{\rm Re}a(P_{{C_H}}f, (I-P_{{C_H}})f)\geq 0\quad \hbox{for all }f\in V_Y.$$ 

% By definition of $V_Y$, one sees that the projection $P_{{C_H}}$ satisfies $P_{{C_H}}V_Y \subset V_Y$
% if and only if
% \begin{enumerate}[(i)]
% \item $P_{{C_H}}H^1(\Omega,W) \subset H^1(\Omega,W),$
% \item $ P_C Y \subset Y.$
% \end{enumerate}
% \footnote{\textcolor{red}{wie robin gemerkt hat, die \"aquivalenz ist eigentlich nicht ganz klar -- aber die implikation von unten nach oben und die implikation von oben nach (ii). was tun? sich damit begn\"ugen? in konkreten f\"alle ist $P_{{C_H}}H^1(\Omega,W) \subset H^1(\Omega,W)$ ja immer erf\"ullt...}
% \textcolor{blue}{proviamo a dimostrare l'implicazione da sotto a sopra. Sia $f\in V_Y$. Allora $f \in H^1$. Quindi per (i) $P_{C_H}f$ \'e anche in $H^1$. 
% D'altra parte, $f_{|\partial \Omega}(z) \in Y$. Allora per (ii) anche $P_C f(z)$ \'e in $Y$. Ma $P_C(f(z)) = P_{C_H} f(z)$ e quindi $P_{C_H} f(z) \in Y$. Ne concludiamo che $P_{C_H}f $ \'e $H^1$ e ha valori al bordo in $Y$. Quindi \'e in $V_Y$.  Per dimostrare da sopra a 2 si prenda una funzione in $\phi \in H^1(\Omega)$ che non sia costantemente 0 sul bordo, 
% cio\'e tale che esiste $z \in \partial \Omega$ e $\phi(z)\neq 0$.
% Per un arbitrario $w \in Y$ si consideri $w\otimes\phi \in V_Y$. Allora $Y \ni P_{C_H} w \otimes \phi(z)=P_C(w\otimes \phi (z))=P_C(\phi (z)w )$.
% Per l'arbitrariet\'a di $w \in Y$ otteniamo l'invarianza di $Y$.}
% }
% These latter conditions can be checked applying Lemma~\ref{verkett} and t
We recall the following result, which has been stated in~\cite{Bar95} and explicitly proved in~\cite[Lemma~2.3]{ManVogVoi05}.

\begin{lemma}\label{projectioninclusion}
Let $C_1$ and $C_2$ be closed convex subsets of a Hilbert space. 
Then the inclusion $P_{C_1} C_2 \subset C_2$ holds if and only if the inclusion $P_{C_2} C_1 \subset C_1$ holds if and only $P_{C_2}P_{C_1}=P_{C_1}P_{C_2}$.
\end{lemma}

Thus, we obtain the following.

\begin{prop}\label{ordersubspaceI}
Let $C_W$ be a closed convex subset of $W$.
Assume that $0\in C_W$ or that $\Omega$ has finite measure.
Then the inclusion $P_{C_H} V_Y \subset V_Y$ holds if and only if the inclusion $P_ Y C_W \subset C_W$ holds.
\end{prop}

\begin{proof}
Since orthogonal projections of a Hilbert space onto a closed convex subset are Lipschitz continuous, by Lemma~\ref{verkett} one has that the pointwise projections do not affect weak differentiability. That is, $P_{C_H}f\in V=H^1(\Omega;W)$ for all $f\in V$.
As a consequence, we obtain that $P_{{C_H}}V_Y\subset V_Y$ if and only if the boundary condition 
$(P_{{C_H}}f)(z)\in Y$ is satisfied for a.e.\ $z\in\partial\Omega$ and for all $f\in V_Y$, i.e., if and only if for all $f\in V$
\begin{equation}\label{quasipoisson2}
f(z)\in Y \hbox{ for a.e.\ }z\in\partial\Omega\quad\hbox{ implies}\quad 
P_C(f(z))\in Y\hbox{ for a.e.\ }z\in\partial\Omega.
\end{equation}
By Lemma~\ref{projectioninclusion} it suffices to show that this is equivalent to $P_C Y\subset Y$. In fact, if $P_C Y\subset Y$ holds, then~\eqref{quasipoisson2} is satisfied.

Conversely, assume that~\eqref{quasipoisson2} holds and let $y\in Y$ and $z \in \partial \Omega$. 
Fix now a function $f$ of class $H^1(\Omega; \mathbb R)$ such that $f(w)=1$ for all $w\in\partial\Omega$ in a neighborhood of $z$. Then the function $g:={f}\otimes y$ is in $V_Y$. Thus $P_{C}(f(z))=P_C y\in Y$, i.e., $P_C Y\subset Y$.
\end{proof}

If in particular $C_W$ is a subspace, we can characterize a symmetry of the problem, i.e.,\ the invariance of $C_H$ under $(e^{tA_Y})_{t\ge 0}$, as follows.

\begin{prop}\label{invsubsp}
Let $C_W$ be a closed subspace of $W$. Then ${C_H}$ is left invariant under $(e^{tA_Y})_{t\ge 0}$ if and only if
\begin{enumerate}[(1)]
\item the inclusion $P_ Y C_W \subset C_W$ holds,
\item the inclusion $D(x) C^n_W \subset C^n_W$ holds for a.e.\ $x\in \Omega$, and
\item the semigroup generated by $S$ leaves invariant $C_{\partial H}$.
\end{enumerate}
\end{prop}

Here $C_W^n$ denotes the Cartesian product of $n$ copies of $C_W$, a subspace of $W^n$.

\begin{proof}
By Proposition~\ref{ordersubspaceI} it suffices to show that ${\rm Re}a(P_{{C_H}}f, (I-P_{{C_H}})f)\geq 0$ if and only if (2)--(3) hold.\\
Since $C$ is a closed subspace, by linearity $\partial_k P_{{C_H}} f=P_{{C_H}}\partial_k f$ for all $k=1,\ldots,N$, as well as $\hat{T}(P_{C_H}f)=P_{C_{\partial H}} \hat{T}$. Thus,
\begin{eqnarray*}
{\rm Re}a(P_{{C_H}}f, (I-P_{{C_H}})f)&=&{\rm Re}\sum_{h,k=1}^n \int_\Omega (D_{hk}(x)P_{C_H}(\partial_k f(x))|(I-P_{C_H})(\partial_h f(x)))_W dx\\
&&\qquad+{\rm Re}\int_{\partial\Omega} (SP_{C_{\partial H}}(f(z))|(I-P_{C_{\partial H}})(f(z)))_W d\sigma(z).
\end{eqnarray*}
Since the partial derivatives are locally independent from the boundary values, and viceversa, the right hand side is positive if and only if both addends are -- in fact, by linearity one sees that positivity already implies that both terms have to vanish. By a localization argument we conclude that
$$\sum_{h,k=1}^n(D_{hk}(x)P_{C}(\partial_k f(x))|(I-P_C)(\partial_h f(x)))_W=0\qquad \hbox{for all }f\in V\hbox{ and a.e.\ }x\in \Omega$$
as well as
$$(SP_{C_H}(f(z))|(I-P_{C_H})(\partial_h f(z)))_W=0\qquad \hbox{for all }f\in V\hbox{and a.e.\ }z\in \partial\Omega$$
if and only if
\begin{equation}\label{orthog}
(D(x)P_{C_W}z|(I-P_{C_W})z)_{W^n}=0\qquad \hbox{for all }z\in W^n\hbox{ and a.e.\ }x\in\Omega 
\end{equation}
as well as
$$(SP_{C_{\partial H}}w|(I-P_{C_{\partial H}})w)_W=0\qquad \hbox{for all }w\in W.$$
Again by Ouhabaz's criterion, the latter condition is equivalent to invariance of ${C_{\partial H}}$ under the (semi)group generated by the bounded operator $S$. Since moreover
$$(D(x)P_{C}z|(I-P_C)z)_W=((I-P_C)D(x)P_{C}z|z)_W,$$
it is also clear that~\eqref{orthog} is equivalent to the inclusion in condition (2) of the statement.\end{proof}

\section{Positivity and irreducibility}\label{posdom}

Throughout this and the next section, we assume $W$ to be a Hilbert lattice. It is known that each Hilbert lattice is isometrically lattice isomorphic to a Lebesgue space $L^2(X)$ for some measure space $X$, see e.g.~\cite[Cor.~2.7.5]{Mey91}. If moreover $W$, hence $L^2(X)$ and also $L^1(X)$ are separable, then one sees that the measure can be taken to $\sigma$-finite (in fact, even finite).

As a consequence, if $f\in H$, $f(x)$ can be regarded as a function $X \ni y \mapsto (f(x))(y) \in \mathbb C$.
For this reason, we sometimes write, with a slight abuse of notation
$$
f(x,y):=(f(x))(y).
$$
In other words, one can regard
$H=L^2(\Omega;W)$ as a scalar valued Lebesgue space 
$$H\cong L^2(\Omega)\otimes L^2(X)\cong L^2(\Omega\times X)$$
on a $\sigma$-finite measure space. Also observe that the lattice structure of $W$ permits to consider
the notion of a \emph{local} bounded linear operator
on $W$\footnote{Recall that by~\cite[Thm.~2.3]{AreTho05} linear bounded local operators on a vector-valued $L^p$-space are exactly multiplication operators by a suitable operator-valued function of class $L^\infty(\Omega;{\mathcal L}_s(W))$.}.

Recall that a vector lattice is an ordered vector space such that $x\vee y:=\sup\{x,y\}$ is well-defined for all vectors $x,y$. Accordingly the notions of positive and negative parts $x^+,x^-$ as well as of absolute value $|x|$ are well-defined, cf.~\cite[\S~C-I]{Nag86}.  A Hilbert lattice is a vector lattice  \emph{and} a  Hilbert space endowed with a lattice norm, i.e., with a norm $\|\cdot\|$ such that $|x|\leq |y|$ implies $\|x\|\leq \|y\|$, where $|\cdot|$ denotes the absolute value on the vector lattice.

We are thus promptly led to introduce similar concepts for vectors of $H$, provided $W$ is a vector lattice. More precisely, let us define the positive and negative parts
$$f^+:\Omega\ni x\mapsto (f(x))^+\in W,\qquad  f^-:\Omega\ni x\mapsto (f(x))^-\in W$$
as well as the absolute value
$$|f|:\Omega\ni x\mapsto |f(x)|\in W.$$
One sees that they all belong to $H$, provided that $f\in H$. In particular, also $H$ is a Hilbert lattice. Since $\|\cdot\|_W$ is a lattice norm, we deduce that also $\|\cdot\|_H$ is a lattice norm, hence $H$ is a Hilbert lattice, % i.e., $|f|\leq |g|$ implies $\|f\|_W\leq \|g\|_W$, 
where of course
\begin{equation}\label{functnorm}
\|f\|_H:\Omega\ni x\mapsto \|f(x)\|_H=\left(\int_X \| f(x,y)\|^2 dy\right)^\frac{1}{2} \in \mathbb R.
\end{equation}

In the following we consider the bounded or unbounded \emph{order intervals}
\begin{eqnarray*}
{[a,b]_W}&:=&\left\{w \in L^2(X): a\le f(x)\le b\hbox{ for a.e.\ }x\in X\right\},\\
{[a,\infty)_W}&:=&\left\{w \in L^2(X): a\le f(x)\hbox{ for a.e.\ x}\in X\right\},\\
{(-\infty,b]_W}&:=&\left\{w \in L^2(X): f(x)\le b\hbox{ for a.e.\ }x\in X\right\},
\end{eqnarray*}
in $W=L^2(X)$, where $a,b\in W$. Of course, they are closed convex subsets of $W$ to which we can apply the invariance theory developed in~\cite[\S~2.1]{Ouh05}.

Having defined order intervals in $W$, it is natural to extend this notion to $H=L^2(\Omega;W)$ by setting
\begin{eqnarray*}
{[a,b]}_H&:=&\left\{f \in H: f(x)\in [a,b]_W\hbox{ for a.e.\ }x\in \Omega\right\},\\
{[a,\infty)}_H&:=&\left\{f \in H: f(x)\in [a,\infty)_W\hbox{ for a.e.\ x}\in \Omega\right\},\\
{(-\infty,b]}_H&:=&\left\{f \in H: f(x)\in (-\infty,b]_W\hbox{ for a.e.\ }x\in \Omega\right\}.
\end{eqnarray*}
We say that $[a,b]_W$ is the order interval \emph{corresponding to ${[a,b]}_H$}, and likewise for the unbounded order intervals. Similarly we can consider the \emph{corresponding} order intervals in $\partial H=L^2(\partial\Omega;W)$.

To begin with, we recover a result analogous to Proposition~\ref{ordersubspaceI}. In the following, we denote by $J_W$, $J_H$, and $J_{\partial H}$ any of the three above kinds of order intervals in $W$, $H$, and $\partial H$, respectively.

% \begin{proof}
% The proof essentially mimicks that of  Proposition~\ref{ordersubspaceI}. We only have to check that $P_{J_H}H^1(\Omega;{\mathbb C})\subset H^1(\Omega;{\mathbb C})$. A careful examination of the proofs of~\cite[Lemma~7.5 and Lemma~7.6]{GilTru01},  shows that the assertions therein also hold in the more general context of arbitrary order intervals of vector-valued Sobolev spaces. In particular, if $C$ is an order interval, then again $P_{J_H}H^1(\Omega;{\mathbb C})\subset H^1(\Omega;{\mathbb C})$ since $\nabla P_{J_H} f={\mathbf 1}_{\{f\in C_H\}}\otimes \nabla f$.
% \end{proof}

Let $a_{Y}$ be defined as in~\eqref{formdomain}--\eqref{aform}. Aim of the following is to discuss the issue of invariance of order intervals of this kind under the semigroup $(e^{tA_Y})_{t\geq 0}$. This is e.g.\ relevant because, of course, characterizing whether the semigroup is positive or contractive with respect to the sup-norm amounts to discussing invariance of ${[0,\infty)}_H$ or of ${[-1,1]}_H$ -- by duality, this also yields contractivity with respect to the $L^1$-norm, which is relevant in applications. We emphasize that for a given order interval $J_W$ and some $f\in H$ the derivatives of the functions $P_{J_H}f,(I-P_{J_H})f\in H$ need not have disjoint support, unlike in the scalar-valued case. This may jeopardize the semigroup's positivity, as we see in the following.

\begin{prop}\label{positive1}
Consider an order interval $J_W$ in $W$ containing $0$ and the corresponding order interval $J_H$ in $H$.
Then the following assertions hold.
\begin{enumerate}[(1)]
\item If the semigroup $(e^{tA_Y})_{t\geq 0}$ leaves invariant $J_H$, then $P_Y$ leaves invariant $J_W$.
\item Assume that $D_{hk}(x)$ is a local operator for all $h,k=1,\ldots,n$ and a.e.\ $x\in\Omega$, then the semigroup
$(e^{tA_Y})_{t\geq 0}$ leaves invariant $J_H$ if and only if the semigroup on $\partial H$ generated by $-S$ and the orthogonal projection $P_Y$ leave invariant $J_{\partial H}$ and $J_W$, respectively.
% \item If $S$ is a local operator, then the semigroup $(e^{tA_Y})_{t\geq 0}$ leaves invariant an order interval of $H$ if and only if $D_{hk}(x)$ is a local operator for all $h,k=1,\ldots,n$ and a.e.\ $x\in\Omega$ and the orthogonal projection $P_Y$ leaves invariant the corresponding order interval in $W$.
%\footnote{\textcolor{red}{in~\cite{CarMugNit08} haben wir bewiesen, dass positivit\"at der halbgruppe impliziert, dass die koeffizientenmatrix diagonal ist -- was im diesen kontext wohl so \"ubersetzt werden k\"onnte, dass die koeffizientenoperatoren $D_{hk}$ notwendigerweise lokal sind, sobald die halbgruppe positiv ist. stimmt das?}}
\end{enumerate}
\end{prop}

\begin{proof}
To begin with,  observe that by Proposition~\ref{ordersubspaceI} for any given order interval $J_W$, the inclusion $P_{J_H} V_Y \subset V_Y$ holds if and only if the inclusion $P_ Y J_W \subset J_W$ holds. Furthermore, we obtain the pointwise identity
$$(\partial_h P_{J_H} f)(x,y)={\mathbf 1}_{\{f(x,y)\in J_W(y)\}}\otimes \partial_h f(x,y), \qquad h=1,\ldots,n,\; x\in\Omega,\; y\in X,$$
where $J_W(y)$ denotes the (real) interval $[a(y),b(y)]$ if $J_W=[a,b]_W$ is a bounded order interval, and analogously in the unbounded case. Here $\partial_h$ denotes the partial derivative in the $h^{\rm th}$ direction with respect to the first coordinate, i.e., $x$. Accordingly, it turns out that $(\partial_k P_{J_H} f),(f-\partial_h P_{J_H} f)$, and due to locality also $(D_{hk}\partial_k P_{J_H} f),(f-\partial_h P_{J_H} f)$, do have disjoint supports as function of two variables over the measure space $\Omega\times X$.

First assume $(e^{tA_Y})_{t\geq 0}$ to leave invariant an order interval $J_H$ of $H$. Then $V_Y$ is invariant under $P_{J_H}$. Our preliminary remark yields (1).

In order to prove (2) assume now that $D_{hk}(x)$ is a local operator for all $h,k=1,\ldots,n$ and a.e.\ $x\in\Omega$. Then
\begin{eqnarray*}
{\rm Re}a(P_{J_H}f,f-P_{J_H}f) & = &{\rm Re}\sum_{h,k=1}^n \int_\Omega (D_{hk}(x)\partial_k (P_{J_H}f)(x)| \partial_h (f-P_{J_H}f)(x))_{W}dx \\
&&\qquad+{\rm Re}\int_{\partial\Omega}((S P_{J_{\partial H}}f)(z)| (f-P_{J_{\partial H}}f)(z))_W d\sigma(z)\\
& = &{\rm Re}\sum_{h,k=1}^n \int_\Omega\int_X (D_{hk}(x)\partial_k P_{J_H}f)(x,y)\overline{\partial_h (f-P_{J_H}f)(x,y)}dx dy \\
&&\qquad+{\rm Re}\int_{\partial\Omega}((S P_{J_{\partial H}}f)(z)| (f-P_{J_{\partial H}}f)(z))_W d\sigma(z)\\
&=&{\rm Re}\int_{\partial\Omega}((S P_{J_{\partial H}}f)(z)| (f-P_{J_{\partial H}}f)(z))_W d\sigma(z),
\end{eqnarray*}
where we are using the fact that the trace operator $\hat{T}$ satisfies $P_{J_{\partial H}} \hat{T}= \hat{T}P_{J_H}$. This shows that the condition ${\rm Re}a(P_{J_H}f,f-P_{J_H}f) \ge 0$ is equivalent to the fact that the semigroup generated by $-S$ on $\partial H$ leaves $J_{\partial H}$ invariant.
% Finally, assume $S$ to be a local operator. Computing as in the proof of (1), one obtains that the condition ${\rm Re}a(P_{J_H}f,f-P_{J_H}f)\ge 0$  is equivalent to
% $${\rm Re}\sum_{h,k=1}^n \int_\Omega (D_{hk}(x)\partial_k (P_{J_H}f)(x)| \partial_h (f-P_{J_H}f)(x))_{W}dx\ge 0$$
% and this is in turn equivalent to the locality of $D$.
% This shows (3)\footnote{\textcolor{red}{'stocazzo!}}.
\end{proof}

\begin{exa}
Let $D=Id$ and $Sf:=sf$ for $s:\partial\Omega\to\mathbb R$. Consider Examples~\ref{orthogonaldirichlet} and~\ref{dihedron}. In the first case, it is well-known that the Laplacian with Robin boundary condition generates a positive and $L^\infty$-contractive semigroup, and it is easy to check that in fact Proposition~\ref{positive1}.(2) applies.

In the latter example, one sees that the projection of $W={\mathbb C}^m$ onto $Y=(1,\ldots,1)$ is the matrix
$$P_Y=\frac{1}{m}\begin{pmatrix}
1 & \ldots & 1\\
\vdots & \ddots & \vdots\\
1 & \ldots & 1
                 \end{pmatrix},$$
which clearly leaves the order interval $[0,\infty)_W$ invariant. Hence, Proposition~\ref{positive1}.(2) yields that Kirchhoff boundary conditions give rise to a positive semigroup. On the other hand, anti-Kirchhoff boundary conditions are \emph{not} associated with a positive semigroup, since the projection 
$$P_{Y^\perp}=I-P_Y=\frac{1}{m}\begin{pmatrix}
m-1 & \ldots & -1\\
\vdots & \ddots & \vdots\\
-1 & \ldots & m-1\end{pmatrix}$$
of $W$ onto $Y^\perp$ is not a positive operator.

Regarding $L^\infty$-contractivity, one sees that a matrix leaves invariant $[-1,1]_{{\mathbb C}^m}$ if and only if the sums of the absolute values of its entries sum up at most to 1 on each row. Thus, $P_Y$ does leave $[-1,1]_{{\mathbb C}^m}$ invariant. However, $P_{Y^\perp}$ leaves $[-1,1]_{{\mathbb C}^m}$ invariant if and only if $m\le 2$. This fully characterizes $L^\infty$-contractivity of the diffusion process on an open book with Kirchhoff and anti-Kirchhoff boundary conditions.
\end{exa}

Observe that since $W$ is assumed to be a Hilbert lattice, and hence in particular to be order-complete, the notion of ${\rm sign }u$ of a vector $u\in W$ is well-defined, cf.~\cite[\S\S~C-I.8 and C-II.2]{Nag86}. Therefore also ${\rm sign }f:\Omega\ni x\mapsto {\rm sign }(f(x))$ is well-defined.

We recall the definition of closed ideal of a Hilbert lattice, cf.~\cite[Def.~2.19]{Ouh05}.

\begin{defi}\label{closed idealdefinition}
Let $X,Y$ be subspaces of a Hilbert lattice. Then $X$ is called an \emph{closed ideal} of $Y$ if
\begin{itemize}
\item $x\in X$ implies $|x|\in Y$ and
\item $x\in X$, $y\in Y$, and $|y|\leq |x|$ imply $y\;{\rm sgn}x\in X$.
\end{itemize}
\end{defi}

\begin{rem}\label{antikk}
Observe that we are not requiring that $X$ is a subspace of $Y$. For example, $\{0\}\times\{0\},\{0\}\times{\mathbb C},{\mathbb C}\times\{0\}$ are all closed ideals of ${\mathbb C}^2$, but also $\langle(-1,1)\rangle$ is a closed ideal of $\langle(1,1)\rangle$.
\end{rem}

We conclude this section briefly discussing irreducibility issues. Recall that an operator on the Hilbert lattice $H$ is called \emph{irreducible} if the only closed ideals of $H$ it leaves invariant are the trivial ones.

In the context of scalar-valued equations, i.e., if the state space is $H=L^2(\Omega;\mathbb C)$,  the closed ideals of $L^2(\Omega;\mathbb C)$ are exactly those closed subspaces of of the form $L^2(\omega,\mathbb C)$, where $\omega \subset \Omega$ is a measurable set. In the vector-valued context the situation is more involved. 
\begin{enumerate}[(1)]
\item First, all closed ideals appearing in the scalar-valued case define in a natural way  closed ideals of $L^2(\Omega;W)$, i.e., for each measurable $\omega \subset \Omega$,
the closed subspace $L^2(\omega, W)$ is a closed ideal of $L^2(\Omega;W)$. 

\item Then, for each closed ideal $I$ of $W$ the subspace 
\begin{equation*}
L^2(\Omega; W) \supset L^2(\Omega; I) := \{f \in L^2(\Omega; W) : f(z) \in I \hbox{ for a.e.\ } z \in \Omega\}
\end{equation*} 
is a closed ideal of $L^2(\Omega;W)$. 
\item Finally, if $\Omega\ni x\mapsto\mathcal P(x)\in L^\infty(\Omega; {\mathcal L}_s(W))$ is a weakly measurable function such that ${\mathcal P}(x)$ is an orthogonal projection onto a closed ideal of $W$ for a.e.\ $x\in \Omega$, then the subspace
\begin{equation}
I_\mathcal P := \{f \in L^2(\Omega;W): f(x) \in {\rm Range \,}\mathcal P(x) \hbox{ for a.e.\ } x \in \Omega \}
\end{equation}
is also a closed ideal of $L^2(\Omega;W)$.
\end{enumerate}

These notions of closed ideal of $L^2(\Omega;W)$ are not equivalent: of course, the closed ideals considered in (1)--(2) are special cases of those considered in (3). In fact, it has been proved in~\cite{CarMug09b} that each closed ideal of $H$ is of the form (3). Such a general notion is not necessarily the most convenient one, though. In fact, much seems to depend on the underlying model.

The notion in (2) seems to be more natural in the context of parabolic systems (and in particular of $1$-dim networks), whereas that in (3) is the proper one in the general context of parabolic equations with infinite dimensional state spaces.

Assume e.g.\ we want to describe diffusion on an open book with $N$ pages: then it seems that the relevant space is $(L^2(\Omega;{\mathbb C}))^N$, implying that the correct notion of closed ideal is that of a functions space over a \emph{subset of  a page of the open book} as in (2) above. In other words, it is reasonable to say that the associated semigroup is irreducible whenever the heat localized inside a single page will be trasmitted to further pages  -- with continuity and Kirchhoff-type conditions, this is the case if and only if the open book is connected, i.e., if and only if all pages are glued to the binding. %Observe that this need not be the case even under the (standing) assumption that $\Omega$ is connected.

More generally, if $I_W$ is a closed ideal of $W$ and we want to discuss invariance of $I_H$, then we can promptly apply Proposition~\ref{invsubsp}. Assume in particular the operators $D(x)$, $x\in\Omega$, as well as $S$ to be multiples of the identity, so that conditions (2) and (3) of Proposition~\ref{invsubsp} are clearly satisfied for all subspaces. Then the diffusion semigroup $(e^{ta_Y})_{t\ge 0}$ is irreducible (in the sense of (2)) if and only if $P_Y I_W \not\subset I_W$ for any closed ideal $I_W$ of $W$, i.e., if and only if the operator $P_Y$ is irreducible. (In the case of an open book with $N$ pages, this of course means that the $(e^{ta_Y})_{t\ge 0}$ is irreducible if and only the $N\times N$-matrix $P_Y$ is not similar to a block upper triangular matrix via a permutation).

\medskip
If we are instead interested in the invariance in the sense of (1), then we promptly observe that no (nontrivial) closed ideal of type (1) is left invariant under the action of any of the semigroups we consider throughout this paper. In fact, the projection onto $L^2(\omega, W)$ is the pointwise multiplication with the characteristic function of $\omega$ and this operator does not leave $H^1(\Omega; W)$ invariant, hence the semigroup $(e^{ta_Y})_{t\ge 0}$ is irreducible for any $Y$, in particular on open books. (Observe that Proposition~\ref{ordersubspaceI} does not apply, since $L^2(\omega;W)$ is not of the form $I_H$ for any closed ideal $I_W$ of $W$).

\medskip
Finally, the general case of irreducibility is more involved. It has been observed in~\cite{CarMug09b} that $(e^{ta_Y})_{t\ge 0}$ is irreducible in the sense of (3) if and only if $H$ is 1-dimensional.

\section{Domination issues}\label{dominff}

The aim of this section is to present several results concerning domination of semigroups associated with the form considered throughout this paper, both by means of semigroups associated with $(a_Y,V_Y)$ for different $Y$ (in Subsection~\ref{domin1}) and by means of a diffusion semigroup acting on a space of scalar-valued functions (in Subsection~\ref{domin2}) .

As in the previous section, let $W$ and hence $L^2(\Omega;W)$ be Hilbert lattices.

\subsection{Domination by vector-valued semigroups}\label{domin1}

The following result is~\cite[Formula~(2.7)]{Ouh05}.

\begin{lemma}\label{dominc}
The set ${\mathcal C}:=\{(x,y)\in W\times W: |x|\leq y\}$ is closed and convex in $W\times W$ and the orthogonal projection onto ${\mathcal C}$ is given by
$$P_{\mathcal C}(x,y):=\frac{1}{2}((|x|+|x|\wedge {\rm Re}\, y)^+ {\rm sgn}x,(|x|\vee {\rm Re}\,y+{\rm Re}\,y)^+).$$
\end{lemma}

\begin{theo}\label{dominvector}
Let $Y_1,Y_2$ two subspaces of $W$ with $Y_1\hookrightarrow Y_2$. Assume $P_{Y_2}$ to be a positive operator. Let finally $\rho_1,\rho_2$ be $L^\infty(\partial\Omega;{\mathcal L}_s(W))$-functions such that $\rho_1(z),\rho_2(z)$ are positive operators for a.e.\ $z\in \partial\Omega$. Consider forms $a_1,a_2$  defined by
$$a_1(f,g)=\int_\Omega \left(D(x)\nabla f(x) | \nabla g(x)\right)_{W^n} dx+\int_{\partial\Omega} \left(\rho_1(z)f(z) | g(z)\right)_W d\sigma(z),\qquad f,g\in V_{Y_1}$$
and
$$a_2(f,g)=\int_\Omega \left(D(x)\nabla f(x) | \nabla g(x)\right)_{W^n} dx+\int_{\partial\Omega} \left(\rho_2(z)f(z) | g(z)\right)_W d\sigma(z),\qquad f,g\in V_{Y_2}.$$
Then the semigroup $(e^{ta_1})_{t\geq 0}$ is dominated by $(e^{ta_2})_{t\geq 0}$, i.e.
$$|e^{ta_1}f(x,y)| \le e^{ta_2}|f|(x,y),\qquad t\ge 0,\; f\in H,\; x\in\Omega,\; y\in X,$$
 if and only if $Y_1$ is a closed ideal of $Y_2$ and moreover $\rho_1(z)\ge \rho_2(z)$ in the sense of positive operators for a.e.\ $z\in\partial\Omega$.
\end{theo}

Observe that under the assumptions of the theorem the semigroup $(e^{ta_2})_{t\ge 0}$ is positive.

\begin{proof}
Consider the subset ${\mathcal C}_H:=\{(f,g)\in H\times H: |f|\leq g\}=\{(f,g)\in H\times H: (f(x),g(x))\in {\mathcal C}\hbox{ for a.e.\ }x\in\Omega\}$. Then ${\mathcal C}_H$ is closed and convex. By Lemma~\ref{dominc} the projection of $H\times H$ onto ${\mathcal C}_H$ is given by
$$P_{{\mathcal C}_H}(f,g):=\frac{1}{2}((|f|+|g|\wedge {\rm Re}\, g)^+ {\rm sgn}f,(|f|\vee {\rm Re}\,g+{\rm Re}\,g)^+).$$
Consider the form $a_W$ with maximal domain $H^1(\Omega;W)$, i.e., the form with (decoupled) Neumann boundary conditions. 
Due to positivity of $(e^{ta_W})$ and Lemma~\ref{dominc}, one can reason as in~\cite[Cor.~2.22]{Ouh05} and conclude that the domination of $(e^{ta_1})_{t\geq 0}$ by $(e^{ta_2})_{t\geq 0}$ is equivalent to $V_{Y_1}$ being a closed ideal of $V_{Y_2}$ and additionally ${\rm Re} a_1(u,v) \ge a_2(|u|,|v|)$ for all $u,v\in V_{Y_1}$ such that $(u|v)_W=|u|_W |v|_W$.\\
First of all, taking into account surjectivity of the trace operator $H^1(\Omega;W)\partial H^1(\partial\Omega;W)$ and density of $H^1(\partial\Omega;W)$ in ${L^2(\partial\Omega;W)}$, a direct computation shows that $V_{Y_1}$ is a closed ideal of $V_{Y_2}$ if and only if $Y_1$ is a closed ideal of $Y_2$; and moreover ${\rm Re} a_1(u,v) \ge a_2(|u|_W,|v|_W)$ for all $u,v\in V_{Y_1}$ such that $(u|v)_W=|u|_W |v|_W$ if and only if $\left(\rho_1(z)x | y\right)_W\ge \left(\rho_2(z)|x|||y|\right)_W$, for all $x,y\in {Y_1}$ such that $(x|y)_W=|x|_W |y|_W$, i.e., if and only if $\rho_1(z)\ge \rho_2(z)$ for a.e.\ $z\in\partial\Omega$.
\end{proof}

\begin{exa}
Fix $k\in \mathbb N$ and let $W=\mathbb C^k$. Consider  the half-space $\Omega:=\{(x_1,\ldots,x_N): x_1 >0\}$ and
 the form $a_{Y_1}$ defined as in Theorem~\ref{dominvector} with $Y_1=\langle \mathsf 1 \rangle$, $D(x)=Id$, $\rho_i(x)=0$, for a.e.\ $x \in \Omega$.
There are no nontrivial subspaces of $W$ $\langle \mathsf 1 \rangle$ is a closed ideal of: on one hand, observe that $\langle \mathsf 1\rangle \subset Y$,  by the first condition in Definition~\ref{closed idealdefinition}. On the other hand, for all $y\in Y$ there exists $x \in \langle \mathsf 1 \rangle$ such that $|y|\leq |x|$.
Since $\langle \mathsf 1\rangle$ is a closed ideal of $Y$, it follows from the second condition in Definition~\ref{closed idealdefinition} that $y \in \langle \mathsf 1\rangle$, i.e.\ $Y=\langle 1\rangle$ as soon as $\langle 1\rangle$ is a closed ideal of $Y$.

As a consequence, no semigroup associated with the form $a_2$ (for any $Y_2$ and $\rho_2$!) defined as in Theorem~\ref{dominvector} dominates the semigroup with boundary conditions defined by $Y_1=\langle \mathsf 1 \rangle$.

One could summarize this by saying that \emph{among the class of boundary conditions on structures of open-book-type, continuity and Kirchhoff boundary conditions define the maximal semigroup}.

Finally, take $k=2$. Observe that interpreting the Kirchhoff boundary conditions as the natural generalization for the Neumann boundary conditions in the context of ramified structures (like open books), an open-book-analog of the conclusion of~\cite{AreWar03} fails to hold: actually, between Dirichlet and continuity + Kirchhoff there are local boundary conditions beyond (decoupled) Robin ones. These are e.g.\ given by so-called anti-Kirchhoff boundary conditions, defined by replacing $\langle 1\rangle$ with $\langle 1\rangle^\perp$, cf.\ Remark~\ref{antikk} -- observe in particular that $\langle 1\rangle^\perp$ is a closed ideal of $\langle 1\rangle$.
\end{exa}

\begin{exa}
In particular, the semigroup with Dirichlet boundary conditions (i.e., associated with $(a_{\{0\}},V_{\{0\}})$ is dominated by any other semigroup of the class considered in the previous section. The semigroup $(e^{ta_Y})_{t\geq 0}$ is dominated by the semigroup governing a system of uncoupled diffusion equations with Neumann boundary conditions, corresponding to $(e^{ta_W})_{t\geq 0}$, if and only if $Y$ is a closed ideal of $W$. This in particular shows that the conclusion of~\cite{AreWar03} carries over to the vector valued case: between Dirichlet and Neumann there is no local boundary condition beyond decoupled Robin ones.
\end{exa}

\subsection{Domination by scalar-valued semigroups}\label{domin2}

We want to discuss domination of $(e^{ta})_{t\ge 0}$ by some semigroup generated by the common Laplace operator on the scalar-valued space $L^2(\Omega;{\mathbb C})$. To this aim, we apply some ideas presented in~\cite[\S~4]{ManVogVoi05}. 

Following~\cite[Def.~3.3]{ManVogVoi05} we introduce the notion of a generalized closed ideal.  In the following we use the notation introduced in~\eqref{functnorm}.

\begin{defi}\label{generalizedide}
A subspace $H_0$ of $L^2(\Omega;W)$ is called a \emph{generalized closed ideal} of a subspace $\mathfrak h_0$ of  $L^2({\Omega};{\mathbb C})$ if
\begin{itemize}
\item $f\in H_0$ implies $\|f\|_W\in \mathfrak h_0$ and
\item $f\in H_0$, $g\in \mathfrak h_0$, and $|g|\leq \|f\|_W$ imply $g\otimes{\rm sign}f\in H_0$.
\end{itemize}
\end{defi}

In particular, $L^2(\Omega;W)$ is a \emph{generalized closed ideal} of $L^2({\Omega};{\mathbb C})$.

The following arguments are inspired by~\cite[\S~3]{ManVogVoi05} and~\cite[\S~2.5]{Ouh05}. 

\begin{theo}\label{dominscalar}
Let $d\in L^\infty(\Omega;M_{n}(\mathbb C))$ satisfy pointwise a uniform ellipticity condition. Let the coefficient $D_{hk}$ be given by $d_{hk}\otimes Id$, i.e., $(D_{hk}f)(x)=d_{hk}(x)f(x)$ for all $h,k=1,\ldots,n$ and a.e.\ $x\in\Omega$. Let moreover $r\in L^\infty(\partial\Omega;\mathbb C)$. Let the coefficient $D_{hk}$ be given by $d_{hk}\otimes Id$, i.e., $(D_{hk}f)(x)=d_{hk}(x)f(x)$ for all $h,k=1,\ldots,n$ and a.e.\ $x\in\Omega$, and similarly $S=s\otimes Id$, i.e., $(Sf)(z)=s(z)f(z)$ for a.e.\ $z\in \partial \Omega$.

Define the forms
$$b(f,g):=\int_{\Omega} (d(x)\nabla f(x)|\nabla g(x))_{\mathbb C^n}dx+\int_{\partial\Omega} r(z) f(z)\overline{g(z)}d\sigma(z),\qquad f,g\in H^1(\Omega;{\mathbb C}),$$
and 

$$a_W(f,g):=\int_{\Omega} ((D\nabla f)(x)|\nabla g(x))_{W^n}dx +\int_{\partial\Omega} ((S f)(z)|{g(z)})_W d\sigma(z),\qquad f,g\in H^1(\Omega;W).$$

Then $(e^{ta_W})_{t\ge 0}$ is dominated by $(e^{-tb})_{t\ge 0}$, i.e.,
$$
\|e^{ta_W}f\|_W \leq e^{-tb}\|f\|_W, \qquad t \geq 0, f \in H.
$$
\end{theo}

\begin{proof}
By~\cite[Thm.~4.1]{ManVogVoi05} it suffices to check that $V=H^1(\Omega;W)$ is a generalized closed ideal of $H^1(\Omega; \mathbb C)$ and that
$a(f,g) \geq b(\|f\|,\|g\|)$ for all $f,g \in V$ such that 
\begin{equation}\label{conditionide}
(f(x) | g(x))_W= \|f(x)\|_W\|g(x)\|_W\qquad\hbox{ for a.e.\ }x\in\Omega. 
\end{equation}

We first observe that by definition $\|f\|_W\in H^1(\Omega;{\mathbb C})$ as soon as $f \in V$: thus, $V$ is a generalized closed ideal of $H^1(\Omega; \mathbb C)$.

In order to check the second condition, we compute
\begin{eqnarray*}
b(\|f\|_W,\|g\|_W)
&= & \int_{\Omega} (d(x)\nabla\|f\|_W(x) | \nabla\|g\|_W(x))_{\mathbb C^n} dx\\ 
&&\qquad +\int_{\partial\Omega}  r(z)\|f(z)\|_W \|g(z)\|_W d\sigma(z)\\ 
& =& \int_{\Omega} ({\rm sign}f(x) | {\rm sign}g(x))_W ((D\nabla f)(x) | \nabla g(x))_{W^n} dx \\
&& \qquad +\int_{\partial\Omega} ({\rm sign}f(z) | {\rm sign}g(z))_W ((Sf)(z) | g(z))_{W} d\sigma(z) \\
&=&\int_\Omega ((D\nabla f)(x) | \nabla g(x))_{W^n} dx +\int_{\partial\Omega} ((Sf)(z) | g(z))_{W} d\sigma(z)\\
&=& a_W(f,g),
\end{eqnarray*}
where the equality in the third line holds by~\eqref{conditionide}.
\end{proof}

\begin{rem}
It is known that if $r\ge0$, then the semigroup associated with $b$ (i.e., the semigroup generated by a general elliptic operator with Robin boundary conditions in the scalar-valued case) satisfies Gaussian estimates under several geometric conditions on $\partial\Omega$, including $\Omega$ being Lipschitz (which is our standing assumption). Accordingly, also all the semigroups dominated by $(e^{-tb})_{t\ge 0}$ satisfy the same estimates (with same constants). In particular, under the assumptions of Theorem~\ref{dominscalar},  $(e^{ta_W})_{t\ge 0}$ satisfies Gaussian estimates, i.e., there exist $c,d>0$ such that
$$
\|e^{ta_W}f\|_W \leq cG(dt)\|f\|_W, \qquad t \in [0,1], f \in H
$$
where $(G(t))_{t \geq 0}$ denotes the Gaussian semigroup on $\mathbb R^N$. Kernel estimates like this one are an important tool in discussing spectral and regularity properties, cf.~\cite{Are04}.

Combining this with Theorem~\ref{dominvector} one concludes that all semigroups associated with a form $(a_Y,V_Y)$ defined as in Theorem~\ref{dominscalar} satisfy Gaussian estimates, provided $Y$ is a closed ideal of $W$. This is of course not surprising, since in this case the system decouples in subsystems of diffusion equations equipped with Dirichlet or Robin boundary conditions, each of which has Gaussian estimates.
\end{rem}

\section{Appendix}
\subsection{Extension of operators}\label{extensapp}
Let $X$ be a $\sigma$-finite measure space and fix a Hilbert basis $(e_n)_{n\in \mathbb N}$ of $W$. 
Then for all $f\in L^2(X;W)$ there exists uniquely determined functions $f_n \in L^2(X)$ such that
\begin{equation}\label{decomposition}
f(x)= \sum_{n \in N} f_n(x) e_n.
\end{equation}
Furthermore,
\begin{equation}\label{decomposition2}
\|f\|^2_{L^2(X,W)}= \sum_{n \in \mathbb N} \|f_n\|^2_{L^2(X)}.
\end{equation}
To see this, define
$$
f_n(x):= ( f(x)| e_n)_W.
$$
Computing
$$
\| f \|^2_{L^2(X,W)}= \int_{X}\|f(x)\|^2_W dx= \int_X \sum_{n \in \mathbb N} |(f(x) | e_n)_W|^2 dx= \sum_{n \in \mathbb N} \int_X |(f(x)|e_n)_W|^2 dx
$$
shows that $f_n \in L^2(X)$ for all $n \in \mathbb N$ and that the decomposition~\eqref{decomposition} holds.

Assume now that there exist functions $g_n$ with the property~\eqref{decomposition}.
Then
$$
\sum_{n \in N} (f_n-g_n)(x) e_n= 0
$$
implies that $f_n=g_n$ a.e.

\begin{theo}\label{extension}
For an open domain $\Omega$ and $s \ge 0$, consider $T \in \mathcal L (H^s(\Omega;{\mathbb C}), L^2(\partial \Omega;{\mathbb C}))$.
Define 
$$
\hat{T}f(x):= \sum_{n \in \mathbb N} Tf_n(x) e_n,\qquad f \in H^s(\Omega;W), x \in \Omega.
$$
Then
$$\hat{T} \in \mathcal L (H^s(\Omega; W),L^2(\partial \Omega; W))$$
and
$$\|T\|=\|\hat{T}\|.$$
\end{theo}

A similar result holds if $L^2(\partial \Omega;W)$ is replaced by $H^r(\Omega;W)$, $r\ge 0$, or more generally for a pair of Sobolev spaces $(H^s(\Gamma_1;{\mathbb C}), H^r(\Gamma_2;{\mathbb C}))$ defined on two sufficiently smooth manifolds $\Gamma_1,\Gamma_2$. We omit the details.

%\textcolor{red}{mi sono reso conto che nella dimostrazione non utilizziamo da nessuna parte propriet\`a geometriche di $\Omega$. la mia sensazione \`e che il principio di estensione valga per qualsiasi coppia di spazi di sobolev generici. azzardo una formulazione che secondo me non ha bisogno di alcun cambiamento nella dimostrazione: se sei d'accordo, rimpiazza il vecchio teorema col nuovo. l'unica questione \`e: che ipotesi di regolarit\`a sulle variet\`a $\Gamma$ ci occorrono per formulare il teorema, introdurre gli spazi di sobolev e per poter considerare $L^2(\Gamma)$ come uno spazio di misura $\sigma$-finito? una variet\`a compatta ovviamente basterebbe, ma \`e un'ipotesi troppo forte. idee?}
%\textcolor{blue}{In Gilbarg-Trudinger c'\'e molto a proposito. Perch\'e non prendiamo semplicemente il setting che scelgono loro?}

\begin{proof}
We first assume $s=0$. By~\eqref{decomposition}
$$
\|\hat{T} f \|^2_{L^2(\partial\Omega; W)}= 
\sum_{n \in \mathbb N} \| T f_n\|^2_{L^2(\partial\Omega;{\mathbb C})} \leq 
\|T\|^2 \|f\|^2_{L^2(\Omega; W)},
$$
for all $f\in L^2(\Omega;W)$, and this shows $\|\hat{T}\|\leq \|T\|$.
For arbitrary $(x_n)_{n \in \mathbb N} \in \ell^2(\mathbb N)$ and $g \in L^2(\Omega;{\mathbb C})$ such that $\|(x_n)_{n \in \mathbb N}\|_{\ell^2(\mathbb N)}=\|g\|_{L^2(\Omega;{\mathbb C})}=1$ define $g_n(x):= x_n g(x)$ and $\hat{g}(x):=\sum_{n \in \mathbb N} g_n(x) e_n$. By~\eqref{decomposition2} we have 
$\| \hat{g} \|^2_{L^2(\Omega; W)}=1$ and therefore
$$
\|\hat{T}\|^2\geq \|\hat{T} \hat{g}\|_{L^2(\partial\Omega;W)}^2 = \sum_{n \in \mathbb N} \|T g_n\|^2_{L^2(\partial\Omega;W)}
= \sum_{n \in \mathbb N} |x_n|^2 \| T g\|^2_{L^2(\Omega;{\mathbb C})}= \|Tg\|^2_{L^2(\Omega;{\mathbb C})}.
$$
The estimate holds for all $g$ such that $\|g\|_{L^2(\Omega;{\mathbb C})}=1$. Thus,
$$
\|\hat{T}\|^2 \geq \sup_{\|g\|=1} \|Tg\|^2_{L^2(\Omega;{\mathbb C})}= \|T\|^2.
$$
Let now $r\in\mathbb N$: then the claim holds by repeating the proofs replacing the function $f$ by its appropriate partial derivatives. The general case of $r\ge 0$ follows by interpolation arguments based on the Riesz--Thorin Theorem.
\end{proof}

As an application of Theorem~\ref{extension} we prove the existence of a vector-valued trace operator.

\begin{exa}\label{tracexa}
Let $\Omega$ be an open Lipschitz domain of ${\mathbb R}^n$. It is classical that there exists a trace operator $T \in \mathcal L(H^1(\Omega;{\mathbb C}), L^2(\partial \Omega;{\mathbb C}))$ 
such that $Tf=f_{|\partial \Omega}$ if $f \in H^1(\Omega;{\mathbb C})\cap C(\overline{\Omega};{\mathbb C})$, see e.g,~\cite[Thm.~1.5.1.3]{Gri85}.
By Theorem~\ref{extension}, this operator can be extended to an operator 
$\hat{T} \in \mathcal L(H^1(\Omega; W), L^2(\partial \Omega; W))$, 
and in fact we denote with an abuse of notation $f_{|\partial \Omega}:=\hat{T}f$ for all $f \in H^1(\Omega;W)$.
\end{exa}

\subsection{Chain rule for vector-valued function}

The following observation is likely to be already known, but we could not find any reference for it. 

\begin{lemma}\label{verkett}
Let $G: W \to W$ be a  Lipschitz continuous mapping and $f \in H^1(\Omega;W)$.
If
\begin{enumerate}
\item $G(0)=0$, or
\item $\Omega$ has bounded measure,
\end{enumerate}
then $G\circ f \in H^1(\Omega;W)$.
\end{lemma}

\begin{proof}
We start observing that the proof of~\cite[Prop.~IX.3]{Bre83} holds also in the vector-valued case with minor changes.
In other words, if $f \in L^2(\Omega;W)$, then $f \in H^1(\Omega;W)$ is equivalent to the existence of a positive constant $C$ with the property that for all open bounded $\omega \subset \Omega$ and all $h\in{\mathbb R}^n$ with $|h| \leq {\rm dist(\omega, \partial \Omega)}$ one has
\begin{equation}\label{BrezisH1}
\int_\omega \|f(x+h)-f(x)\|^2_W dx \leq C|h|^2.
\end{equation}
Assume $G$ to be Lipschitz with constant $L$. First we note that the estimate
$$
\int_\omega \|G(f(x+h))-G(f(x))\|^2_W dx  \leq  \int_\omega L^2\|f(x+h)-f(x)\|^2_W dx
\leq C L^2 |h|^2
$$
holds for every $f \in H^1(\Omega; W)$, since by assumption $f$ satisfies~\eqref{BrezisH1}.
It remains to show that $G\circ f \in L^2(\Omega; W)$. 

If $G(0)=0$, estimate
$$
\int_\Omega \|G(f(x))\|^2_W dx =\int_\Omega \|G(f(x)) - 0\|^2_W dx=\int_\Omega \|G(f(x)) - G(0)\|^2_W dx \leq L^2\int_\Omega \|f(x)\|^2_W dx < \infty.
$$
Thus $G \circ f \in L^2(\Omega;W)$ and the above criterion applies.

Let now $\Omega$ have finite measure. Fix an arbitrary vector $w \in W$ and estimate
\begin{eqnarray*}
\|G(f(x))\|_W & = & \|G(f(x))-G(w) + G(w)\|_W\\
& \leq &  L\|f(x)-w\|_W + \|G(w)\|_W\\
& \leq & L\|f(x)\|_W + L\|w\|_W + \|G(w)\|_W.
\end{eqnarray*}
Squaring and integrating with respect to $x$ we obtain
\begin{eqnarray*}
\int_\Omega \|G(f(x))\|^2_W dx &\leq&
\int_\Omega 
L^2\|f(x) \|^2_W + L^2\|w \|^2_W + \|G(w)\|^2_W
+2L^2\|f(x) \|_W \| w \|_W   \\
&& 
+2L  \|f(x) \|_W\|G(w)\|_W 
+2L  \|w \|_W\|G(w)\|_W 
dx. 
\end{eqnarray*}
The latter integral is finite since $\Omega$ has finite measure and this completes the proof.
\end{proof}

\bibliographystyle{plain} 
\bibliography{literatur.bib}

\end{document}